\documentclass[11pt]{article}
\usepackage{latexsym, amsfonts, amsmath, amsthm, graphicx}                              

\def\~{\tilde }
\def\RR{{\mathbb R}}
\def\II{{\bf 1}}

\theoremstyle{plain}
\newtheorem{theorem}{Theorem}[section]
\newtheorem{lemma}[theorem]{Lemma}
\newtheorem{corollary}[theorem]{Corollary}

\theoremstyle{definition}
\newtheorem{remark}[theorem]{Remark}
\newtheorem{example}[theorem]{Example}

\begin{document}

\title{A sharp isoperimetric bound for convex bodies}

\author {Ravi Montenegro \thanks{School of Mathematics, Georgia Institute of Technology, 
         Atlanta, GA 30332-0160, monteneg@math.gatech.edu; supported in part by
         a VIGRE grant.}}

\date{}

\maketitle

\begin{abstract}
\noindent We consider the problem of lower bounding a generalized Minkowski measure of
subsets of a convex body with a log-concave probability measure, conditioned on the set size. 
A bound is given in terms of diameter and set size, which is sharp for all set sizes, dimensions,
and norms. In the case of uniform density a stronger theorem is shown which is also sharp. 

\vspace{2ex}
\end{abstract}

\section{Introduction} \label{sec:intro}

It is a classic result that among all surfaces in $\RR^3$ enclosing a fixed volume, the
sphere has minimal surface area, as measured by the Minkowski measure $\mu^+$. A related extremal problem 
shows that half spaces minimize surface area for a Gaussian distribution in $\RR^n$ \cite{Bor75.1}. 

One variation on these results is to consider log-concave measures $\mu$ supported on a convex body $K$, 
i.e. a closed and bounded convex set. Recall that the Minkowski measure
$\mu^+(S) = \lim_{h\rightarrow 0^+} \mu(S_h\setminus S)/h$, where $S_h$ denotes
the set of points at most distance $h$ from $S$. 
In certain applications it is better to work in a different norm. We define the generalized
Minkowski measure $\mu^+(S)$ as before, but where distance in $S_h$ is measured in
the desired norm $\|\cdot\|$. If $\|\cdot\|$ is the $\ell_2$ Euclidean norm then this is
the standard Minkowski measure. 

Our main result is the following.

\begin{theorem} \label{thm:isoperimetry}
Let $\mu$ be a log-concave probability measure supported on a convex body $K\subset\RR^n$. 
For all measurable sets $S\subset K$ with $\mu(S)\leq 1/2$ it follows that
$$
(diam\,K)\,\mu^+(S) \geq \mu(S)\,G(1/\mu(S))\ ,
$$
where $(diam\,K)$ is measured in some arbitrary norm $\|\cdot\|$, $\mu^+(S)$ is
the generalized Minkowski measure in that same norm, and $G(1/\mu(S))$ is given
by
\begin{equation} \label{eqn:G}
G(1/x) = \frac{\gamma^2 e^{\gamma}}{e^{\gamma} (\gamma-1)+1}\ ,
\end{equation}
where $\gamma > 0$ is the unique solution to
$$
x = \frac{e^{\gamma} (\gamma-1)+1}{(e^{\gamma}-1)^2} \ .
$$
Moreover, this bound is sharp for every value of $x$, diameter of $K$,
dimension $n$ and norm $\|\cdot\|$.
\end{theorem}

This extends theorems of Dyer and Frieze \cite{DF91.1} and Lov\'asz and Simonovits \cite{LS93.1} 
in which they did not condition on set size. Our sharpness also
shows that their theorems are tight only when $x=1/2$.

The main tool in our proof is the Localization Lemma of Lov\'asz and Simonovits, which 
makes it possible to reduce an $n$-dimensional integration problem into a one dimensional
problem. A unique aspect of our method is that we start with an unknown lower bound, given 
by $G(1/x)$, proceed to discover which properties $G(1/x)$ must have to apply Localization, and only
at the final step, after reducing this to a one dimensional problem, do we determine the function
$G(1/x)$. All other applications of Localization of which we are aware begin with a conjectured 
lower bound and proceed to show it to be correct. However, by not making any assumptions to begin
with we are able to obtain the a sharp lower bound which would have been a very unlikely initial
candidate. 

It does not appear possible to write the function $G(1/x)$ in closed form. However,
in Corollary \ref{cor:logbounds} we show that $G(1/x)$ behaves like $\log(1/x)$,
being bounded below by $2+\log(1/2x)$ and above by $2\,\log_2(1/x)$. It is interesting to 
note that for graphs with a nice geometric structure, as with the grid $[k]^n$ (see \cite{BL91.1}), 
this shows that the graph number and the edge-isoperimetry are likely to
differ by a logarithmic factor.

We are also able to apply our methods to the more specific case of the uniform distribution.
In Theorem \ref{thm:uniform} we give a bound which is again sharp for every set size $x$,
dimension $n$ and norm $\|\cdot\|$. The main improvement is when $x$ is small, and in 
Corollary \ref{cor:uniform} it is shown that when $x<2^{-n}$ then $G(1/x)$
behaves like $n/\sqrt[n]{x}$. Example \ref{ex:hypercube} shows that, at least for the 
$\ell_{\infty}$ norm, the extremal cases on the hypercube $[0,1]^n$ are always within a factor $3$
of the extremal cases on general convex bodies.

In fact, in general, the extremal cases are relatively simple to state.
When $\mu$ is uniform and the enclosed volume $\mu(S) = x$, 
dimension $n$ and diameter $D$ are fixed, then there is a truncated cone $K$ with a subset $S$ that is extremal
(the slope of the cone depends on the dimension $n$).
More generally, when $\mu$ is a log-concave probability measure and $x=\mu(S)$, dimension $n$ and diameter 
$D$ are fixed, then we show 
that the long thin cyclinder $[0,1]\times[0,\epsilon]^{n-1}$ with a one dimensional 
exponential distribution $F(x)=e^{\gamma\,x_1}/(\epsilon^{n-1}\,\int_0^1 e^{\gamma\,x}\,dx)$ contains a subset 
$S=[0,s]\times[0,\epsilon]^{n-1}$ that is extremal as $\epsilon\rightarrow 0^+$, 
where $\gamma$ is from Theorem \ref{thm:isoperimetry}
and both $s$ and $\gamma$ are independent of the dimension $n$.

Since Theorem \ref{thm:isoperimetry} is sharp, then all bounds of the form 
$$
(diam\,K)\,\mu^+(S) \geq f(\mu(S))
$$
will follow as corollaries. For instance, 
$$
(diam\,K)\,\mu^+(S) \geq
\mu(S)\,\mu(K\setminus S)
 \ \left(4 + \log\left(1/4\,\mu(S)\,\mu(K\setminus S)\right)\right)  
$$
strengthens a result of Kannan, Lov\'asz and Montenegro \cite{KLM03.1}.
A different weakening leads to a Gaussian isoperimetric function,
\begin{eqnarray} \label{eqn:gaussian}
(diam\,K)\,\mu^+(S) \geq \sqrt{2\pi}\,I_{\gamma}(\mu(S))\ ,
\end{eqnarray}
where $\mu(S) \leq 1/2$ and $I_{\gamma}(x)$ is the Gaussian 
isoperimetric function (see (\ref{eqn:gaussian2})).

We note that other authors
\cite{Bob00.1,KLS99.1} have proven related results in which quantities measuring the well-roundedness of 
$K$ were fixed, rather than the diameter $D$, but their results appear to be tight only in asymptotics
and not when conditioned on dimensions and volumes $\mu(S)$ as is the case here.
For instance, Bobkov \cite{Bob00.1} used a Pr\'ekopa-Leindler inequality to obtain a related result.
\begin{theorem} \label{thm:Prekopa}
Let $\mu$ be a log-concave probability measure in $\RR^n$. For all measurable sets $S\subset \RR^n$, for every point 
$x_0\in \RR^n$, for every number $r > 0$, and for standard $\ell_2$ Minkowski measure,
$$
2r\,\mu^+(S) \geq \ \mu(S)\,\log\frac{1}{\mu(S)} + (1-\mu(S))\,\log\frac{1}{1-\mu(S)} + \log \mu\{|x-x_0|\leq r\}\ .
$$
\end{theorem} 

This is not directly comparable with our result because Theorem \ref{thm:Prekopa} considers shape (via $r$) 
as well as set size $\mu(S)$. However, if $r$ is a radius then the log term drops out.
In this case a comparison can be made and our result is stronger. Of course, $r$ can be chosen so that the 
$\log$ term need not drop out completely, in which case the bounds are not comparable.

The paper proceeds as follows. In Section \ref{sec:inequality} we prove 
Theorem \ref{thm:isoperimetry}. Section \ref{sec:bounds} proves various bounds
on the quantity $G(1/x)$. We conclude with the uniform case in section \ref{sec:uniform}.


\section{The Proof} \label{sec:inequality}

Recall that a function $f : \RR^n \rightarrow \RR^+$ is {\em log-concave} if 
$\forall x,\,y\in \RR^n,\,t\in[0,1] : f[tx + (1-t)y] \geq f(x)^t\,f(y)^{1-t}$,
i.e., $\log f$ is a concave function on on the support of $f$. 
In particular, non-negative concave functions are log-concave.

A measure $\mu$ is log-concave if for every measurable
$A,\,B\subset \RR^n :\, \mu(t\,A+(1-t)\,B) \geq \mu(A)^t\,\mu(B)^{1-t}$.
All log-concave measures are induced by log-concave functions, so that 
$\mu$ is log-concave if and only if there is a log-concave function $F$ such that
for every measurable $S\subset \RR^n:\,\mu(S) = \mu_F(S) = \int_S F(x)\, dx$.


A {\em lower semi-continuous} function is one which is a limit of a monotone
increasing sequence of continuous functions. For example, the indicator of an open
set, or the negative of the indicator of a closed set.

The lemma below is a variation on results in \cite{KLS99.1,LS93.1}. 
\begin{lemma}[Localization Lemma] \label{lem:localization}
Let $g$ and $h$ be lower semi-continuous Lebesgue integrable functions on
$\RR^n$ such that
$$
\int_{\RR^n} g(x)\,dx \geq 0
\quad and \quad
\int_{\RR^n} h(x)\,dx = 0\ .
$$
Then there exist two points $a,\,b\in \RR^n$ and a linear function 
$\ell : [0,1]\rightarrow \RR^+$ such that
$$
\int_0^1 \ell(t)^{n-1}\,g((1-t)a + tb)\,dt \geq 0
\quad and \quad
\int_0^1 \ell(t)^{n-1}\,h((1-t)a + tb)\,dt = 0\ .
$$
\end{lemma}

We begin by reducing our problem into a one-dimensional one. 

\begin{theorem} \label{thm:conditions}
Let $\mu_F$ be a log-concave measure on $\RR^n$ induced by log-concave function $F$,
with compact support $K$ and a disjoint partition $K= S_1\cup S_2\cup B$ with 
$\mu_F(S_1) \leq \mu_F(S_2)$. Also, let $t$ and $d$ be such that 
$d \geq diam\,K$, $t \leq dist(S_1,S_2)$, both relative to some norm $\|\cdot\|$.

If $G : [2,\,\infty) \rightarrow \RR^+$ and $x\,G(1/x) : (0,1/2]\rightarrow \RR^+$ 
are monotonically non-decreasing, then
\begin{equation} \label{eqn:isoineq_logconcave}
\frac {d-t}{t}\,\mu_F(B) \geq \mu_F(S_1)\,
       G\left(\frac{\mu_F(K\setminus B)}{\mu_F(S_1)}\right)
\end{equation}
holds for all such partitions 
if it holds for all one-dimensional probability distributions 
$\~F(t) = e^{\gamma\,t}$ and all intervals 
$S_1=[0,s),\,B=[s,s+t],\,S_2=(s+t,1]$ such that $\mu_{\~F}(S_1) \leq \mu_{\~F}(S_2)$.
\end{theorem}

\begin{remark}
The dependence on $B$ in the lower bound 
    can be removed by replacing $d-t$ with $d$ and $K\setminus B$ by $K$.
    However, this leads to a slightly weaker result which no longer 
    extends \cite{LS93.1} in terms of $t$.
\end{remark}

\begin{proof}
Assume a contradiction, i.e. $\exists F,\,K,\, S_1,\, S_2,\, B$ with
\begin{equation} \label{eqn:counterexample}
\frac {d-t}{t}\,\frac{\mu_F(B)}{\mu_F(K\setminus B)} 
   < \frac{\mu_F(S_1)}{\mu_F(K\setminus B)}\,G\left(\frac{\mu_F(K\setminus B)}{\mu_F(S_1)}\right)\ .
\end{equation}
By continuity of measure,
if $S_1$ and $S_2$ are increased by a small enough amount, with 
$B=K\setminus \left(S_1\cup S_2\right)$ and $t$ decreasing accordingly, then this
still gives a counterexample. It can then be assumed that $S_1$ and $S_2$ are
open, with $B$ closed.

The Localization Lemma can be used to reduce this to a one-dimensional problem. 
The following two conditions 
will decrease the left side of the counterexample, while keeping the right side constant.
\begin{eqnarray*}
g(t) &=& F(t)\,\left(A\,\II_{K\setminus B}(t) - \II_{B}(t)\right) \textrm{ where } A = \mu(B)/\mu(K\setminus B) \\
\textrm{ and } h(t) &=& F(t)\,\left(x\,\II_{K\setminus B}(t) - \II_{\overline{S_1}}(t)\right) \textrm{ where } x = \mu(S_1)/\mu(K\setminus B)\ ,
\end{eqnarray*}
where $\overline{S_1}$ indicates the closure of $S_1$.
The condition $\int g(t)\,dt\geq 0$ assures that $\mu(B) / \mu(K\setminus B)$ not increase when changing to one-dimension,
while $\int h(t)\,dt=0$ causes $x = \mu(S_1) / \mu(K\setminus B)$ to stay constant, and these two conditions
imply $\mu(S_1)\leq\mu(S_2)$ in the one-dimensional problem as well. 
This decreases the left side of 
(\ref{eqn:counterexample}) while increasing the right side (by the conditions on $G$), and hence gives a
one-dimensional counterexample. The one dimensional problem has smaller diameter (length)
than $K$ and larger separation $t$, so the same $d$ and $t$ are valid in the one-dimensional problem.  
Moreover, by linearity all norms are equivalent in $\RR^1$ up to a constant factor; these constants
cancel out when taking $(d-t)/t$, so it can be assumed that the norm is standard Euclidean length. 
Without loss, assume the one dimensional problem is on $[0,1]$. 

In the sequel, for $t\in[0,1]$ then $F(t)$ will denote the restriction to the one-dimensional problem,
i.e. $\ell(t)^{n-1}\,F(t\,a+(1-t)\,b)$, and for $[u,v]\subseteq [0,1]$ then
$\mu_F([u,v])=\int_u^v \ell(t)^{n-1}\,F(ta+(1-t)b)\,dt$.

Suppose there is a one-dimensional counterexample (with $<$ in (\ref{eqn:isoineq_logconcave})) 
where $B$ consists of a single interval, i.e. $[0,1]=[0,u)\cup[u,v]\cup(v,1]$
where $S_1=[0,u)$, $B=[u,v]$, $S_2=(v,1]$. Let $\log \~F(t)$ be the line
$\log \~F(t) = A + \gamma\,t$ passing through the points $(u,\,\log F(u))$ and
$(v,\, \log F(v))$. Log-concavity of $F(t)$ implies that $\~F(t) \leq F(t)$ in $B$
and $\~F(t) \geq F(t)$ in $S_1 \cup S_2$. Then $\mu_{\~F}(B) \leq \mu_{F}(B)$ and the left side of the counterexample
decreases in changing from $\mu_F(B)$ to $\mu_{\~F}(B)$. 
Also, $\mu_{\~F}(S_1) \geq \mu_{F}(S_1)$ and $\mu_{\~F}(S_2) \geq \mu_{F}(S_2) \geq \mu_F(S_1)$, 
so that by monotonicity of $G(x)$ the right side of the counterexample increases in going from $\mu_F(K\setminus B)$ to 
$\mu_{\~F}(K\setminus B)$, and by monotonicity of $x\,G(1/x)$ there is another increase in going from
$\mu_F(S_1)$ to $\min\{\mu_{\~F}(S_1),\,\mu_{\~F}(S_2)\}$.  This completes the single interval case.

In general, the one-dimensional problem may have many intervals. We use a trick of
Lov\'asz and Simonovits \cite{LS93.1} to reduce the general case to the single interval case.
Suppose that $\mu_F([0,r)) > \mu_F((s,1])$ for the leftmost maximal interval
$[r,s]\subseteq B\subset [0,1]$, or 
$\mu_F([0,r)) \leq \mu_F((s,1])$ for the rightmost maximal interval. In the former
case the result follows from the single interval case applied to $S_2=[0,r),\ B=[r,s],\ S_1=(s,1]$,
while the latter case is similar. 
Otherwise, consider consecutive maximal intervals $[r,s]$ and $[u,v]$ of $B$ such that
$\mu_F([0,r))\leq\mu_F((s, 1])$ but $\mu_F([0,u))>\mu_F((v,1])$. 
Either $S_1$ or $S_2$ is a subset of $[0,r)\cup (v, 1]$, assume $S_1$.
If the single interval case has been proven then
\begin{eqnarray*}
\frac {1-t}{t}\,\mu_F([r,s]) &\geq&
   \mu_F([0,r))\,G\left(\frac{\mu_F(K\setminus[r,s])}{\mu_F([0,r))}\right) \\
&\geq&
   \mu_F([0,r)\cap S_1)\, G\left(\frac{\mu_F(K\setminus B)}{\mu_F([0,r)\cap S_1)}\right)  \\
&\geq&
   \mu_F([0,r)\cap S_1)\, G\left(\frac{\mu_F(K\setminus B)}{\mu_F(S_1)}\right)
\end{eqnarray*}
where the inequalities follow from the monotonicity conditions on $G(x)$ and $xG(1/x)$. Likewise,
$$
\frac {1-t}{t}\,\mu_F([u,v]) \geq
    \mu_F((v,1]\cap S_1)\, G\left(\frac{\mu_F(K\setminus B)}{\mu_F(S_1)}\right)\ .
$$

Adding these expressions together gives
\begin{eqnarray*}
\frac {1-t}{t}\,\mu_F(B) &\geq& \frac{1-t}{t}\,\left(\mu_F([r,s]) + \mu_F([u,v])\right) \\
                  &\geq& \mu_F(S_1)\, G\left(\frac{\mu_F(K\setminus B)}{\mu_F(S_1)}\right)\ ,
\end{eqnarray*}
as desired. If it were $S_2 \subseteq [0,r)\cup (v,1]$ then the same steps would hold with
$S_2$. Since $\mu_F(S_2) \geq \mu_F(S_1)$ then the monotonicity of $x\,G(1/x)$ 
would then imply the result for $\mu_F(S_1)$.
\end{proof}

This reduces the problem to a one dimensional one on $[0,1]$, with 
log-concave measure $\mu_F(A) = \int_A e^{\gamma\,t}\,dt$ 
and intervals $S_1 = [0,s),\,B=[s,s+t],\,S_2=(s+t,1]$. We now find the optimal $G(x)$ 
for the one-dimensional problem, which leads to the optimal function 
$G(x)$ in Theorem \ref{thm:conditions} as well.

\begin{theorem} \label{thm:isoperimetry_soln}
Let $G : [2,\infty) \rightarrow \RR$ be defined by $G(2) = 2$ and 
\begin{equation} \label{eqn:gamma_logconcave}
\forall \gamma > 0\,:\
G(1/x) = \frac{\gamma^2 e^{\gamma}}{e^{\gamma} (\gamma-1)+1}
\quad where \quad x = \frac{e^{\gamma} (\gamma-1)+1}{(e^{\gamma}-1)^2} \in (0,1/2) \ .
\end{equation}
Then the conditions of Theorem \ref{thm:conditions} are satisfied,
and the theorem is sharp at every value of $x=\mu(S_1)/\mu(K\setminus B)$.
\end{theorem}

\begin{remark}
The general case of sharpness was given in the introduction.
Recall that $K=[0,1]\times [0,\epsilon]^{n-1}$. When $x=1/2$ then $\gamma=0$ and $F=1$,
with sharpness when $S_1$ is half the cylinder, i.e. $S_1=[0,1/2]\times[0,\epsilon]^{n-1}$.
This is the same as the sharpness result for Dyer and Frieze's \cite{DF91.1} version of 
Theorem \ref{thm:isoperimetry} that does not condition on $x$.
Similarly, when $x=1/2$ and $t > 0$ in Theorem \ref{thm:conditions} then $F=1$ and
$S_1=[0,(1-t)/2)\times [0,1]^{n-1}$, $B=[(1-t)/2,(1+t)/2]\times[0,1]^{n-1}$.
This was the sharp case for Lov\'asz and Simonovits \cite{LS93.1}. Our
bounds equal theirs when $x=1/2$ and are strictly better when $x<1/2$.
\end{remark}

\begin{remark}
Alternatively, $\gamma$ can be interpreted as the slope of $x\,G(1/x)$ because
\begin{eqnarray*} 
\frac{d}{dx} [x\,G(1/x)] 
&=& \frac{d}{d\gamma} \left[ \frac{\gamma^2 e^{\gamma}}{(e^{\gamma}-1)^2}\right]
  \ \left[\frac{dx}{d\gamma}\right]^{-1} \\
&=& 
\frac
{\frac{(2\gamma e^{\gamma} + \gamma^2 e^{\gamma})(e^{\gamma}-1)^2 - \gamma^2 e^{\gamma} 2 e^{\gamma} (e^{\gamma}-1)} {(e^{\gamma}-1)^4} }
{\frac{(e^{\gamma}+e^{\gamma}(\gamma-1))(e^{\gamma}-1)^2 - (e^{\gamma}(\gamma-1)+1) 2 e^{\gamma} (e^{\gamma}-1)} {(e^{\gamma}-1)^4} 
} \\
&=&
\gamma\ .
\end{eqnarray*}
\end{remark}

\begin{proof}[Proof of Theorem \ref{thm:isoperimetry_soln}]
Consider a one-dimensional (single interval) counterexample to Theorem \ref{thm:conditions}
with $\mu_{\~F}(S_1) \leq \mu_{\~F}(S_2)$. The case $\mu_{\~F}(S_2) \leq \mu_{\~F}(S_1)$ follows similarly.
\begin{equation} \label{eqn:integral_counter}
\frac{1-t}{t}\,\frac{\mu_{\~F}(B)}{\mu_{\~F}(K\setminus B)}
< \frac{\mu_{\~F}(S_1)}{\mu_{\~F}(K\setminus B)}\, 
     G\left(\frac{\mu_{\~F}(K\setminus B)}{\mu_{\~F}(S_1)} \right)\ .
\end{equation}

Write the intervals as $S_1 = [0,s)$, $B = [s,s+t]$ and $S_2 = (s+t,1]$. If $t$ is decreased,
while fixing $x = \mu_{\~F}(S_1)/\mu_{\~F}(K\setminus B)$ and adjusting $s$ in order to keep $x$ constant,
then the right side of the counterexample remains constant in $t$.

For the left side, some simple algebra shows that 
$\mu_{\~F}(S_1) = \frac{x}{1-x}\,\mu_{\~F}(S_2)$,
and therefore 
$$
e^{\gamma\,s}-1 = \frac{x}{1-x}\,\left(e^{\gamma} - e^{\gamma\,(s+t)}\right)\ .
$$
Solving for $e^{\gamma\,s}$ proves that
$$
\gamma\,\mu_{\~F}(S_1) = e^{\gamma\,s}-1 = x\,\frac{e^{\gamma}-e^{\gamma\,t}}{1+x\,(e^{\gamma\,t}-1)}\ ,
$$
from which it follows that the left side of the counterexample is 
$$
\frac{1-t}{t}\,\frac{\mu_{\~F}(B)}{\mu_{\~F}(K\setminus B)}
 = \left(1 + x(e^{\gamma}-1)\right)\,\frac{1-t}{t}\,\frac{e^{\gamma\,t}-1}{e^{\gamma}-e^{\gamma\,t}}\ .
$$

By Lemma \ref{lem:min} with $D=e^{\gamma}$ it follows that, taking 
$t \rightarrow 0^+$ on both sides of (\ref{eqn:integral_counter}) 
gives another counterexample:
\begin{equation} \label{eqn:counter}
\gamma x + \frac{\gamma}{e^{\gamma}-1} < x\, G(1/x)\ .
\end{equation}

Fix $x$ and minimize the left side with respect to $\gamma$.
$$
\frac{\partial}{\partial \gamma}
    \left(\gamma x + \frac{\gamma}{e^{\gamma}-1} \right)
= x + \frac{e^{\gamma}-1-\gamma e^{\gamma}}{(e^{\gamma}-1)^2}
$$
When $\gamma > 0$ then this is increasing in $\gamma$ and so the
root is an absolute minimum of (\ref{eqn:counter}), i.e. the 
minimum occurs at the solution to
\begin{equation}  \label{eqn:bijection}
x = \frac{e^{\gamma} (\gamma-1) + 1}{(e^{\gamma}-1)^2} \in (0,\,1/2)\ .
\end{equation}

Observe that when $\gamma \in (0,\infty)$ then (\ref{eqn:bijection}) is a bijection onto
$x\in (0, 1/2)$. Since the solution to (\ref{eqn:bijection}) is the minimum of the 
left side in (\ref{eqn:counter}) then there is another counterexample with
$$
\gamma \frac{e^{\gamma} (\gamma-1) + 1}{(e^{\gamma}-1)^2}
  + \frac{\gamma}{e^{\gamma}-1}
< x\,G\left(\frac{1}{x}\right)\ .
$$
This simplifies to give a contradiction to (\ref{eqn:gamma_logconcave}).

When $\gamma < 0$ then reverse orientation and consider $[1,0]$. This reduces it to the
problem just considered.

The above work shows that for fixed $x$ then when $t \rightarrow 0^+$ the
$\gamma$ given by (\ref{eqn:bijection}) leads to an equality in 
(\ref{eqn:isoineq_logconcave}), so $G$ is sharp.
\end{proof}

\begin{lemma} \label{lem:min}
If $D > 1$ and $t > 0$ then
$$
\frac{1-t}{t}\,\frac{D^t-1}{D-D^t} \geq \frac{\log D}{D-1}\ ,
$$
with the minimum occuring as $t\rightarrow 0^+$.
\end{lemma}

\begin{proof}
%
Let $D = e^{\gamma}$ for some $\gamma > 0$. Cross multiplying and simplifying, it suffices to show
$$
(1-t)\,(e^{\gamma\,t}-1)\,(e^{\gamma}-1) - \gamma t\,(e^{\gamma}-e^{\gamma\,t}) \geq 0\ .
$$

Plugging in the Taylor series for $e^x$ into the left side, and factoring out
$t$ or $1-t$ factors gives
\begin{eqnarray*}
LHS &=& \gamma^2\,t(1-t)\,\sum_{k=0}^{\infty} \frac{(\gamma\,t)^k}{(k+1)!}\,
                  \sum_{k=0}^{\infty} \frac{\gamma^k}{(k+1)!}
  - \gamma^2\,t\,\sum_{k=0}^{\infty} \frac{\gamma^k}{(k+1)!}\,(1-t^{k+1}) \\
 &=& \gamma^2\,t(1-t)\,\left[
     \sum_{k=0}^{\infty} \gamma^k\,\sum_{i=0}^k \frac{t^i}{(i+1)!\,(k-i+1)!} 
   - \sum_{k=0}^{\infty} \frac{\gamma^k}{(k+1)!}\,\sum_{i=0}^k t^i \right] \\
 &=& \gamma^2\,t(1-t)\,\sum_{k=0}^{\infty} \gamma^k\,\sum_{i=0}^k t^i 
         \left[\frac{\binom{k+2}{i+1}}{(k+2)!} - \frac{1}{(k+1)!}\right] 
 \geq 0,
\end{eqnarray*}
where the inequality uses that $\binom{k+2}{i+1} \geq k+2$ for
$i\in\{0,\ldots,\, k\}$.
\end{proof}

\section{Bounding $G(1/x)$} \label{sec:bounds}

Theorem \ref{thm:isoperimetry} gives an optimal bound, but it seems impossible to write $G(1/x)$
in closed form. We give here a few upper and lower bounds which show that $G(x)$ is 
essentially logarithmic in $x$.

\begin{corollary} \label{cor:logbounds}
If $x = \mu(S) \leq 1/2$, then 
\begin{eqnarray*}
4x(1-x)\,\log_2\left(\frac{1}{x(1-x)}\right) 
   \geq& x\,G(1/x) 
   &\geq x(1-x)\,\left(4 + \log\left(\frac{1}{4\,x(1-x)}\right)\right) \\
2x\,\log_2(1/x)
  \geq& x\,G(1/x)
  &\geq x\,(2+\log(1/2x)) 
\end{eqnarray*}
and has limit
$$
\frac{G(1/x)}{\log(1/x)} \xrightarrow{x\rightarrow 0^+} 1 \ .
$$
\end{corollary}


The first lower bound is a stronger form of a result of 
Kannan, Lov\'asz and Montenegro \cite{KLM03.1}.
Computer plots show that the absolute error is no more than $0.0051$, 
or at most $0.51\%$ of the $[0,1]$ range of $\mu^+(S)$,
and the relative error is no more than $7\%$.

Another lower bound of interest is
\begin{equation} \label{eqn:gauss}
(diam\,K)\,\mu^+(S) \geq \sqrt{2\pi}\,I_{\gamma}(x)\ ,
\end{equation}
where $I_{\gamma}(x)$ is the so-called Gaussian isoperimetric function
\begin{equation} \label{eqn:gaussian2}
I_{\gamma}(x) = \varphi\circ \Phi^{-1}(x) \quad where \quad \varphi(t) = \frac{1}{\sqrt{2\pi}}\,e^{-t^2/2}
 \quad and \quad \Phi(t) = \frac{1}{\sqrt{2\pi}}\,\int_{-\infty}^t e^{-y^2/2}\,dy\ ,
\end{equation}
which makes its appearance in many isoperimetric results, such as Bobkov's \cite{Bob00.1}.
This lower bound is weaker than the first one in the corollary and
so we do not prove it here.

\begin{proof}[Proof of Corollary]
The second upper bound follows
because the general log-concave bound is at least as small as the uniform one (see the next
section). By Corollary \ref{cor:uniform}, 
taking $n\rightarrow\infty$, the upper bound on the uniform case for $x>2^{-n}$ becomes a bound 
on the case $x>0$. The first upper bound follows from the second one
because $2x(1-x)\geq x$ and $x(1-x) \leq x$.

For the limiting case
$$
\lim_{x\rightarrow 0^+} \frac{G(1/x)}{\log(1/x)} 
  = \lim_{\gamma\rightarrow\infty} \frac{\frac{\gamma^2 e^{\gamma}}{e^{\gamma}(\gamma-1)+1}}
       {\ln\left(\frac{(e^{\gamma}-1)^2}{e^{\gamma} (\gamma-1) + 1}\right)}
  = 1\ .
$$

To prove the first lower bound, substitute the expression for $\mu^+(S)$ in Theorem \ref{thm:isoperimetry} 
into the lower bound and rearrange terms. This reduces the problem to one of showing
\begin{equation} \label{eqn:KLM}
\frac{\gamma^2\,e^{\gamma}}{(e^{\gamma}-1)^2}\,\frac{1}{x(1-x)} 
- \left(4 + \log\left(\frac{1}{4\,x(1-x)}\right)\right) \geq 0
\quad
where
\quad
x=\frac{e^{\gamma}\,(\gamma-1)+1}{(e^{\gamma}-1)^2}\ .
\end{equation}
In order to show that (\ref{eqn:KLM}) is non-negative, 
it suffices to show that $\frac{d}{d\gamma}\,(Eqn.\,\ref{eqn:KLM}) \geq 0$,
or equivalently that the sign of the derivative is never negative. 
Multiplying the derivative by positive functions does not affect its sign, 
so positive factors can be cancelled out of the derivative before checking its sign.
This differentiation and cancellation of terms can be performed
repeatedly; if the final expression is non-negative, and if after each intermediate derivative 
the value at $0$ was non-negative, then it follows that (\ref{eqn:KLM}) holds.

Consider
$$
\frac{d^4}{d\gamma^4}
  \left[e^{-\gamma}\,\frac{d^2}{d\gamma^2} \left[ e^{-\gamma}\,
  \frac{d^4}{d\gamma^4} \left[ e^{-\gamma} \,
  \frac{d}{d\gamma} \left[ 
  \frac{(e^{\gamma}(\gamma-1)+1)^2\,(e^{\gamma}-(\gamma+1))^2\,(e^{\gamma}-1)}
       {e^{\gamma}(\gamma-2)+(\gamma+2)}
  \, \frac{d}{d\gamma}\,(Eqn.\ref{eqn:KLM})
 \right] \right] \right] \right]
 = 1296\,e^{\gamma}\ .
$$
It can be verified that each intermediate derivative was non-negative as $\gamma\rightarrow 0^+$,
and $1296\,e^{\gamma}$ is trivially non-negative, so by the earlier remarks (\ref{eqn:KLM}) follows.

The second lower bound follows from this by checking that the difference of the two lower bounds
is concave with minima of $0$ at $x\rightarrow 0^+$ and $x=0.5$.
\end{proof}

\section{The Uniform Distribution} \label{sec:uniform}

When the distribution $F$ is uniform over $K$ then the results from the previous sections
can be strengthened slightly. The proof is similar, but without the
reduction from log-concavity to $e^{\gamma\,t}$. Instead, the extremal cases will be 
truncated cones $\{ {\bf x} : \| <x_2,x_3\,\ldots,x_n>\| \leq 1+\gamma\,x_1\}$, which leads
to a more tedious computation.

\begin{theorem} \label{thm:uniform}
Theorem \ref{thm:isoperimetry} holds for the uniform distribution, but
with optimal $G(1/x)$ in dimension $1$ given by $x\,G(1/x)=1$, and in
dimension $n>1$ given by $G(2) = 2$ and
\begin{equation*}
\forall \gamma > 0 :
x\,G(1/x) = \frac{\gamma\,n}{(1+\gamma)^n-1}
  \left[\frac{(1+\gamma)^{n-1} \gamma (n-1)}{(1+\gamma)^{n-1}-1}\right]^{1-1/n} 
\end{equation*}
where
$$
x = \frac{(1+\gamma)^{n-1} \left[\gamma (n-1) -1\right] + 1}
         {\left[(1+\gamma)^{n-1}-1\right]\left[(1+\gamma)^n-1\right]}\in (0,1/2)\ .
$$
\end{theorem}

How much of an improvement does this give over the log-concave result? 
By fixing a constant $\hat\gamma >0$ and setting $\gamma = \hat\gamma/n$,
then as $n\rightarrow\infty$ the bound in Theorem \ref{thm:uniform} converges
to that of the dimension free Theorem \ref{thm:isoperimetry}, 
just with $\hat\gamma$ in place of $\gamma$.

For finite $n$ the main difference is for small values of $x$.
In particular, the limiting cases in Corollaries \ref{cor:logbounds} and \ref{cor:uniform} (see below)
reveal that when $n$ is fixed and $x\rightarrow 0^+$ then $G(1/x)$ for the log-concave case is 
infinitely smaller than for the uniform case. Therefore the log-concave bound is not a good 
approximation of the uniform result on small subsets, although the following corollary does show
that it is a good approximation when $x>2^{-n}$.

\begin{corollary} \label{cor:uniform}
In Theorem \ref{thm:uniform} the quantity $G(1/x)$ is bounded by
$$
\begin{array}{rccclcl}
\displaystyle
\vspace{1ex}
\frac{n}{\sqrt[n]{x}} &\geq& G(1/x) &\geq& \displaystyle \frac 12\,\frac{n}{\sqrt[n]{x}}  && when\ x\leq 2^{-n}, \\
\displaystyle
2\,\log_2(1/x) &\geq& G(1/x) &\geq& 2+\log(1/2x) &\qquad& when\ x> 2^{-n}\ ,
\end{array}
$$
and has limit
$$
\frac{G(1/x)}{n/\sqrt[n]{x}} \xrightarrow{x\rightarrow 0^+} 1 \ .
$$
\end{corollary}

\begin{proof}
For the first upper bound it suffices to give an example satifying the upper bound for every $x$ and $n$.
Consider an $n$-dimensional unit hypercube $[0,1]^n$.
Then, in $\ell_{\infty}$ norm, a subcube embedded in the corner of volume $x$
will have $\mu^+(S)=n\,x^{1-1/n}$, $\mu(S)=x$ and $diam_{\infty}\,K=1$. Therefore,
$G(1/x)\leq (diam\,K)\,\mu^+(S)/\mu(S) = n/\sqrt[n]{x}$.

When $n=1$ then the lower bound is trivial.

When $n>1$ then 
\begin{eqnarray*}
x\,G(1/x) &\geq& \frac{\gamma\,n}{(1+\gamma)^n-1}
  \left[x\,\left((1+\gamma)^n-1\right)\right]^{1-1/n}  \\
&\geq& \frac n2\,x^{1-1/n}\quad if\ \gamma\geq 1\ .
\end{eqnarray*}
The formula for $x$ is monotone decreasing in $\gamma$, so this
implies the lower bound when 
$$
x \leq x_{\gamma=1}\ where\ x_{\gamma=1} = \frac{2^{n-1}\,(n-2)+1}{(2^{n-1}-1)(2^n-1)} \\
 > \frac{2^{n-1}\,(n-2)}{2^{n-1}\,2^n} = \frac{n-2}{2^n} \ .
$$

Then $x_{\gamma=1}>1/2^n$ for $n\geq 3$, and when $n=2$ then $x_{\gamma=1}=1/3 > 1/2^n$ again.

The second lower bound follows from Section \ref{sec:bounds}, as the lower bound
for the uniform problem is certainly no worse than that for the general log-concave.

The second upper bound follows by an example, again.
Once again consider the unit hypercube $[0,1]^n$ with $\ell_{\infty}$ norm, but this time for
$x\in \left(2^{-(k+1)},2^{-k}\right]$ (where $1\leq k<n$) consider the $k$-dimensional
subcube, i.e., $S=[0,x^{1/k}]^k\,\times\,[0,1]^{n-k}$.
Then $\mu(S)=x$, $\mu^+(S)=k\,x^{1-1/k}$ and therefore
$G(1/x) \leq (diam\,K)\,\mu^+(S)/\mu(S) = k/\sqrt[k]{x}$.
The function $\sqrt[k]{x}\,\log_2(1/x)$ is minimized in 
$\left(2^{-(k+1)},2^{-k}\right]$ at $x=2^{-k}$, with minimum
$k/2$. This implies that $G(1/x) \leq 2\,\log_2(1/x)$.

For the limiting case,
\begin{eqnarray*}
\lim_{x \rightarrow 0^+} \frac{G(1/x)}{n /\sqrt[n]{x}}
&=&
\lim_{\gamma\rightarrow\infty}
\frac{\gamma n \sqrt[n]{(1+\gamma)^{n-1}-1}
  \frac{\left[ (1+\gamma)^{n-1} \gamma (n-1)\right]^{1-1/n}}
       {(1+\gamma)^{n-1} \left[\gamma (n-1) - 1\right] + 1}
}
{n \sqrt[n]{\frac
                   {\left[(1+\gamma)^{n-1}-1\right]\left[(1+\gamma)^n-1\right]}
                   {(1+\gamma)^{n-1} \left[\gamma (n-1) -1\right] + 1}
                }
     }
\\
&=&
\lim_{\gamma\rightarrow\infty}
\frac{\gamma}{\sqrt[n]{(1+\gamma)^n-1}}
\left\{
\frac {(1+\gamma)^{n-1} \gamma (n-1)}
      {(1+\gamma)^{n-1} \left[\gamma (n-1) - 1\right] + 1}
\right\}^{1-1/n} \\
&=& 1\ .
\end{eqnarray*}
\end{proof}

%

It is illustrative to compare these bounds to something that is known exactly.

\begin{example} \label{ex:hypercube}
Consider the $n$-dimensional hypercube $[0,1]^n$ with uniform distribution $F=1$, 
$\ell_{\infty}$ 
norm and $S\subset K$ required to have faces parallel to the surfaces (i.e. $\|u\|_1 = 1$). 
Bollob\'as and Leader \cite{BL91.1} studied surfaces of minimal surface area for this problem 
and found that the extremal sets for fixed $x$ are just the $k$-dimensional subcubes that
we used to determine the upper bounds in Corollary \ref{cor:uniform}, i.e.,
\begin{equation}\label{eqn:BL}
(diam\,K)\,\mu^+(S)/\mu(S) \geq \min_{k\in\{1,\ldots,n\}} k/\sqrt[k]{x}\ .
\end{equation}
Simple calculus shows that $k/\sqrt[k]{x}\log(1/x)\geq e$, with
the minimum occuring at $x=e^{-k}$. Therefore, 
$(diam\,K)\,\mu^+(S)/\mu(S) \geq e\,\log(1/x)$, which shows that the best logarithmic
approximation to (\ref{eqn:BL}) is only a factor $e$ larger than the general lower bound of
Corollary \ref{cor:uniform}. Likewise, the upper bound in the corollary is an upper bound 
to (\ref{eqn:BL}) because it was found by fixing $k$ over certain intervals.
\end{example}

Bollob\'as and Leader solved the hypercube problem in order to find an edge-isoperimetric
inequality on the grid $[k]^n$. The bounds of Corollary \ref{cor:uniform} show
that in graphs with a nice geometric structure, such as $[k]^n$, then the graph number
(or cutset expansion) and the edge-isoperimetry are likely to differ by a logarithmic factor
of inverse set size.

\section{Remarks}

The diameter is often a poor measure of the size of a convex body. For instance, the diameter of
$[0,1]^n$ in the standard Euclidean $\ell_2$ norm is $\sqrt{n}$, whereas the average distance of a 
point from the center is much smaller.
It would be nice if the methods of this paper could be used to allow conditioning on set sizes
in results using other measures of diameter, such as the theorem of
Kannan, Lov\'asz and Simonovits \cite{KLS99.1}
$$
\left(\int |x-x_0|\,\mu(dx)\right)\,\mu^+(S) \geq (\log 2)\,\mu(S)\,\mu(K\setminus S) \ .
$$
When $\mu$ is a probability measure then $\int |x-x_0|\,\mu(dx)$ measures the average
radius of the convex body (or probability distribution) centered at $x_0$, and replaces 
the diameter in this paper. A sharp result conditioned on set size would read something like
$$
\left(\int |x-x_0|\,\mu(dx)\right)\,\mu^+(S) \geq \mu(S)\,F(1/\mu(S))\ .
$$
However, a problem arises because
the average radius, $\int |x-x_0|\,\mu(dx)$, may increase
when the reduction is made to a one-dimensional problem. This contrasts to the diameter,
which is non-increasing. Therefore it is necessary to ``waste'' an inequality
in the Localization Lemma to hold down this average radius,
while two more inequalities would be needed to find $F(1/\mu(S))$. Since Localization
only allows for two inequalities then our method fails here.


\bibliographystyle{plain}
\bibliography{../references}

\end{document}